\documentclass[11pt,a4paper]{article}

\usepackage{amsmath}
\usepackage{amssymb}
\usepackage{amsthm}
\usepackage{epsf}
\usepackage{graphicx}
\usepackage{color}
\usepackage{pstricks}
\usepackage{pst-plot}

\usepackage{epsfig}

\usepackage{color}
\usepackage[colorlinks]{hyperref}

\usepackage{algorithmic}
\usepackage{algorithm}
\definecolor{darkred}{rgb}{0.5,0,0}
\definecolor{darkgreen}{rgb}{0,0.5,0}
\definecolor{darkblue}{rgb}{0,0,0.5}
\hypersetup{colorlinks
,linkcolor=darkred
,filecolor=darkgreen
,urlcolor=darkred
,citecolor=darkblue}

\newcommand{\Nset}{\mathbb N}

\theoremstyle{plain}
\newtheorem{theorem}{Theorem}

\theoremstyle{definition}
\newtheorem{definition}{Definition}

\theoremstyle{remark}

\newtheorem{claim}{Claim}

\newcommand{\ignore}[1]{}

\DeclareMathOperator{\lengthofword}{lg}
\DeclareMathOperator{\aprogress}{AP}
\DeclareMathOperator{\encmap}{enc}
\DeclareMathOperator{\topofstack}{first}

\begin{document}
\setlength{\footnotesep}{1em}

\title{\Large Using the Incompressibility Method to obtain Local Lemma results for Ramsey-type Problems}

\author{%
  Pascal Schweitzer\thanks{Max-Planck-Institut f\"ur Informatik, Campus E1 4, 66123 Saarbr\"ucken, Germany. }
}

\date{April 4, 2008
}


\setlength{\parindent}{0cm}

\maketitle

\begin{abstract}
\noindent We reveal a connection between the incompressibility method and the Lov\'{a}sz local lemma in the context of Ramsey theory.
We obtain bounds by repeatedly encoding objects of interest and thereby compressing strings. The method is demonstrated on the example of van der Waerden numbers. It applies to lower bounds of Ramsey numbers, large transitive subtournaments and other Ramsey phenomena as well.\\

\noindent  \emph{Keywords:} Incompressibility; Lov\'{a}sz Local Lemma; Ramsey Theory; Data Compression; Algorithms

\end{abstract}

\section{Introduction}

The incompressibility method is a general tool appearing in proofs in various fields of mathematics. It is based on the incompressibility of strings:  For strings over any finite alphabet, there is no injective map from the strings of length $\ell$ to the strings of length $ \ell' < \ell $. A more sophisticated approach to the matter uses Kolmogorov complexity of strings, to establish that not only the set of strings as whole but also most individual strings turn out to be incompressible. Li and Vit\'{a}nyi  \cite{LiVitanyi} extensively treat the incompressibility method in their introduction to Kolmogorov complexity. It is well understood how results obtained by the probabilistic first moment method \cite{alon_spencer_92} can be understood in the framework of the incompressibility method. In this paper we show how Ramsey results in the fashion of those obtained with the Lov\'{a}sz local lemma can be obtained via the incompressibility method. On the example of van der Waerden numbers $w(k,c)$ we demonstrate how repeated compression yields bounds that depend only on the number of objects that intersect a given one, showing that $w(k;c) \geq  \frac{c^{k-1}}{4\cdot k} \cdot \frac{k-1}{k}$.
The obtained bound is not the best one known. Szab{\'o} \cite{Szabo90} obtained a lower bound (for the case c = 2) by analyzing explicitly how progressions can intersect in addition to an application of the local lemma. He proved that for any $\epsilon >0$ the two color van der Waerden number fullfills $w(k;2) \geq 2^k/k^\epsilon$ for sufficiently large $k$. For specific $k$ even more is known, Berlekamp \cite{Berle} showed that $w(p+1,2)\geq p \cdot 2^p$ for any prime $p$. The method explained here can be adopted straightforwardly to obtain other Ramsey-theoretic lower bounds for example for Ramsey numbers, large transitive subtournaments and the like. See the book by Graham, Rothschild and Spencer \cite{GRS90} for a broad treatment of Ramsey theory. It is not the aim of this paper to produce new bounds for
these problems. Instead, we aim at a better understanding of one of the most prominent tools of the probabilistic method,  the Lov\'{a}sz local lemma, in the framework of the incompressibility method.


\section{Classical application of the incompressibility\\ method}
We briefly describe a typical application of the incompressibility method to obtain a lower bound on van der Waerden numbers that matches the bound obtained by the first moment method, the original basic probabilistic method. Let $s$ be a string of length $n$ over an alphabet $\Sigma$. In the context of Ramsey theory we think of $s$ as a coloring $\chi \colon \{1 \ldots n \} \rightarrow \Sigma$ of integers. We say $s$ is colored with $c$ colors, where $c=|\Sigma|$ is the size of the alphabet. A \emph{$k$-term arithmetic progression with gap $d$} is a set of the form $\{a,a+d,\ldots,a+(k-1)\cdot d\}$ for $a,d\in\Nset$. It is \emph{monochromatic} if the restriction of $\chi$ to it is constant.

Let $w(k;c)$ be the $c$-color van der Waerden number for $k$-term arithmetic progressions. That is, $w(k;c)$ is the smallest positive integer such that any string of length $w(k;c)$ over an alphabet of size $c$ contains a monochromatic arithmetic progressions of length $k$. Van der Waerden's theorem \cite{Waerden} ensures that this number always exists. 

Now, let $s$ be an arbitrary string of length $n=w(k;c)$. This means $s$ contains a monochromatic arithmetic k-term progression. Number all arithmetic progressions of length $k$. There are at most $n^2/k$ such progressions. Encode one such monochromatic progression by $\lceil {\log_c(n^2/k)}\rceil$ characters. Delete it from the string and write the encoding plus the color of the progression with 1 additional character in front of the string. This mapping is injective, therefore, by the incompressibility of strings, the resulting string is not shorter than the one we started with. We get:
\[\lceil \log_c(n^2/k) \rceil+1 + n - k\geq  n\]
which quickly yields: 
\[ w(k;c) = n > \sqrt{k} \cdot c^{\frac{k}{2}-1} \]


\section{Repeated application of the incompressibility \\method}
{
We will now describe how to obtain a bound that matches (up to a factor of $4/e$) the one obtained via Lov\'{a}sz' local lemma.
In this section we will exclusively deal with monochromatic $k$-term arithmetic progressions. We will therefore casually omit the length of the progressions and the fact that it is monochromatic whenever it is evident from the context. For simplicity we will restrict ourselves to the case of two colors $c=2$. We want to strengthen the result from the previous section about van der Waerden numbers by repeated application of the compression. The increased compression ratio is based on the following idea:

Suppose $s$ is a string of length $n=w(k;2)$ that contains only one progression $\aprogress$. We now replace $\aprogress$ with random bits. By this we mean that every single bit of the progression is replaced by a new bit taken uniformly at random. To encode the first progression we need, as before, $\lceil \log(\frac{n^2}{k})\rceil+1$ bits. After we have replaced the progression, we know that a new progression $\aprogress'$ must have been formed within the string. $\aprogress'$ intersects the old progression $\aprogress$ at least in one position. By intersecting we mean here that the arithmetic progressions as a set of positions intersect. This information enables us to encode $\aprogress'$ with less bits. For example we can obtain an encoding with three pieces of information: 
If we know at what terms of $\aprogress$ and $\aprogress'$ they intersect and how big the spread of $\aprogress'$ is, we can determine $\aprogress'$ from $\aprogress$. For these parameters there are $k$,$k$ and $n/(k-1)$  choices respectively. We then need one additional bit to encode the color of the progression.

Figure \ref{fig:intersect} shows an example of two intersecting progressions. 

\begin{figure}[h]
\begin{center}
 \includegraphics[height=2.0cm]{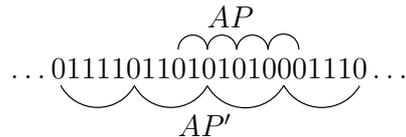}
 \caption{Two 5-term progressions with gaps 2 and 5 intersecting in term 2 of $\aprogress$ and term 3 of $\aprogress'$.}\label{fig:intersect}
 \end{center}
\end{figure}

With this encoding we need all in all $\lceil \log_2(n\cdot k\cdot \frac{k}{k-1}) \rceil +1$ bits to encode $\aprogress'$ together with its color. This is a far less than we needed to encode the first progression $\aprogress$. Note that $\aprogress$ and $\aprogress'$ might be monochromatic at different times during the replacement.

We will now repeatedly use this trick to get arbitrarily close to a compression that encodes $k$ bits with $\lceil \log_2(n\cdot k\cdot \frac{k}{k-1})\rceil$ new bits. For this we first need two formal definitions.

We are only concerned with progressions that are contained within the first $n=w(k;2)$ bits of any string we currently work on. Let $\aprogress$ be such a progression. There are at most $k\cdot k \cdot \frac{n}{k-1}$ progressions that intersect $\aprogress$.
We can encode such a progression $\aprogress'$ by $C = \lceil \log_2(k\cdot k \cdot \frac{n}{k-1}) \rceil$ bits. This encoding then allows us to determine $\aprogress'$ given that we know $\aprogress$.

\begin{definition}[encoding a progression $\aprogress'$ with respect to $\aprogress$]
Fix from now on a mapping $\encmap \colon \mathcal{AP}\sqcap  \mathcal{AP}\rightarrow \{0,1\}^C$, from the set of pairs of intersecting arithmetic progressions within the first $n$ bits to the strings of length $C = \lceil \log_2(k\cdot k \cdot \frac{n}{k-1})\rceil$ which is injective in the first component. The string $\encmap(\aprogress' ,\aprogress)$ is called the \emph{encoding of $\aprogress'$ with respect to $\aprogress$}.
\end{definition}
\begin{definition}[delete/replace a progression]
Let $s =a_1\ldots a_{\lengthofword(s)}$ be a string and $\aprogress = \{ a_t,\ldots,a_{t+k\cdot d}\}$ a progression that is contained in the first $n$ bits. \emph{Replacing $\aprogress$ in $s$} gives a new string $s'$ that is obtained by deleting all elements of $\aprogress$ and then inserting new bits from the end of $s$. More formally we get $s'= \hat{a}_1\ldots \hat{a}_{\lengthofword(s)-k}$ where 
\[\hat{a}_i= \begin{cases}%
\hfill a_{\lengthofword(s)-m} &\text{if } i = t + m \cdot d \text{ for } 0\leq m\leq k-1, \\
\hfill a_i & \text{otherwise.}\\
\end{cases}\]
\end{definition}

We now describe the compression algorithm. Let $s$ be a very long string that is to be compressed. The string $s$ should be longer that $D\cdot k$, where $D$ is the number of replacements that are to be performed. Let $u$ be a string of length $w(k;2)-1$ that contains no arithmetic progression of length $k$. This string exists by the definition of the van der Waerden number. We now work on the string $us$. We will later see that we can compensate for the additional bits from the string $u$ by increasing the number of replacements~$D$. Let $\aprogress_0$ be an arithmetic progression of $us$ among the first $w(k;2)$ bits. This progression will certainly contain the bit at position $w(k;2)$.

Algorithm~\ref{lovasz_algo} is a compression algorithm that takes the string $s$, given on an input tape, and outputs a new string $ps'$, which contains a preamble $p$ that encodes information on how to reobtain the string $s$ from the string $s'$. It manages a queue of progressions, initially only containing $\aprogress_0$, such that at any given time all monochromatic progressions intersect some progression in the queue. The algorithm repeatedly replaces monochromatic progressions that intersect the progression at the head of the queue, and documents all operations as a preamble on the output tape, so that they can be undone. After $D$ replacements, it attaches the remaining input to the back of the output tape. (See Algorithm~\ref{lovasz_algo}.)

\begin{algorithm}
\caption{Repeated compression algorithm}
\label{lovasz_algo}
\begin{algorithmic}[1]
\STATE change input string $s$ to $us$ 
\STATE insert  $\aprogress_0$ to the end of the queue

\WHILE{less than $D$ progressions have been deleted}

\IF {there is a monochromatic progression $\aprogress$ that intersects the progression $\aprogress_{\topofstack{}}$ at the front of the queue } 
	\STATE let $\aprogress$ be the smallest progression (with regard to some fixed ordering) that intersects $\aprogress_{\topofstack{}}$
	\STATE write $\encmap(\aprogress,\aprogress_{\topofstack{}})$
	\STATE write the color of the progression 
	\STATE insert progression $\aprogress$ at the end of the queue
	\STATE delete progression $\aprogress$ on input tape, replace with new bits from the end of the input tape
        \STATE write 1
	
\ELSE
	\STATE remove the first element of the queue
        \STATE write 0
\ENDIF 
\ENDWHILE
\STATE for any progression remaining in the queue write 0 onto the output tape
\STATE copy remaining input to the output tape
\end{algorithmic}
\end{algorithm}

\paragraph{Analysis of algorithm~\ref{lovasz_algo}:}
\paragraph{correctness:}

We do not claim that the output string has any meaning besides encoding the original string. In the view of this, there are only two conditions that might prevent the algorithm from running correctly. First, the replacement might shorten the string beyond a length of $w(k;2)$. This will not happen since we have required the length of $s$ to exceed $D\cdot k$. Second, the queue might be completely emptied and afterwards the first element of the queue might not be well defined. 
 \begin{claim}
 At any time during the execution any monochromatic progression $\aprogress'$ within the first $n$ bits is intersected by some progression $\aprogress$ in the queue.
\end{claim}
As there is always some progression among the first $n$ bits by the claim the queue is never empty. 
We now prove the claim. 
There are two cases why some progression $\aprogress$ might be in the first $n$ bits. Either it has been there from the beginning, it then intersects $\aprogress_0$ in the $n$-th bit, or it developed during the replacement step of some progression $\aprogress'$. But then $\aprogress$ and  $\aprogress'$ intersect, and $\aprogress'$ is enqueued at this point. Since in both cases the respective progression on the stack cannot be removed until $\aprogress$ has been replaced the claim follows.

\paragraph{injectivity:}

We now argue that the mapping that transforms $s$ to the output is injective.
Assume two strings $s$ and $s'$ yield the same output. We can prove by induction that the executions on $s$ and $s'$ at any step yield the same queue as well as the same IF decision (in line 4). The fact that this holds for the IF decision is guaranteed by (and the exact reason why) the algorithm writes 0 and 1 (line 10, 13 and 16) onto the tape depending on the existence of progressions.

Now either $s$ and $s'$ differ in a bit $i$ that is within a progression that is replaced at some point, in which case the progression or the color of it is encoded differently, or it is never deleted at all, in which case it remains within the string and is copied to the output tape in the last step. In both cases th difference on the output tape yields the desired contradiction. 

A formal proof, which we omit here, can be obtained if one defines the space of states as tuples consiting of: A line number that indicates in which step the execution is at some given time, the sequence of progressions that reside on the stack, the content of the input tape and the content of the output tape. One then has to prove that every execution of a single line is reversible, which shows the injectivity.

\paragraph{compression factor:}

We will now compare the length of the output string to the length of the input string:

Say that we have replaced $D$ progressions. This means we have deleted $D \cdot k$ bits from the string. We have written $D$ times 1 (line 10) for the case that we found a progression, since every 1 is followed by a deletion of a progression from the original string. We have written  $D$ times 0  since any progression in the queue leads to exactly one 0 being written, either if it is removed (line 13) or if it remains in the end (line 16).
The map $\encmap()$ encodes each progressions with $ \lceil \log_2 (k\cdot k \cdot \frac{n}{k-1})\rceil$ bits. Additionally we use one bit to encode the color of a progression. Therefore $ \lceil \log_2 (k\cdot k \cdot \frac{n}{k-1})\rceil +1$ bits are written per deleted progression.
In total the new string has length: $\lengthofword(s) + (w(k;2)-1) +  D\cdot (\lceil \log_2 (k\cdot k \cdot \frac{n}{k-1})\rceil+3 ) - D \cdot k $.

Observe that length of the output string depends only on the length of the input string.

The difference to the old string is therefore  $(w(k;2)-1) + D\cdot (\lceil \log_2 (k\cdot k \cdot \frac{n}{k-1})\rceil+3 ) - D \cdot k  $ which must be, by the incompressibility of strings, at least 0, in other words:

\[(w(k;2)-1) + D\cdot (\lceil \log_2 (k\cdot k \cdot \frac{n}{k-1})\rceil+3 ) \geq D \cdot k . \]
Since this holds for arbitrary D we get:
\[\lceil \log_2 (k\cdot k \cdot \frac{n}{k-1})\rceil +3 \geq k \,, \]
which then gives:
\[n > \frac{2^{k-4}}{k} \cdot \frac{k-1}{k} \, .\]

Having analyzed the algorithm we now turn to an improvement of it. First of all the fact that the number of bits needed has to be rounded to an integer is merely an artifact of encoding the possibilities as string. If we reinterpret the inequalities as a form of double counting possibilities it will vanish. In our case there is second more direct way to see this. Instead of repeatedly encoding progressions we can encode the entirety of these progressions at once. This will then yield the inequality 

\[(w(k;2)-1) + \lceil  D\cdot (\log_2 (k\cdot k \cdot \frac{n}{k-1})+3 )\rceil \geq D \cdot k \,. \]

It strengthens the bound by a factor of two. Secondly the method can also be used to bound multicolored van der Waerden numbers from below. Analyzing all steps in the previous section with $c$ colors in mind instead of two, we get that:

\begin{theorem}
The van der Waerdern number $w(k;c)$ for $c$-colored $k$-term arithmetic progressions is bounded below by

\[w(k;c) \geq  \frac{c^{k-3}}{k} \cdot \frac{k-1}{k} = \frac{c^{k-3}}{k} (1-o(1)).\]


\end{theorem}

The obtained bound matches the one originally obtained by Lov\'{a}sz and is worse than later applications of the local lemma by a factor of~$e/4$.
Note that as mentioned in the introduction even these are not the best bounds known. Consult \cite{landman2004} for recent Ramsey theory on the integers.

The local lemma and the method described in this section have in common that they reduce to intersecting objects, since coloring disjoint objects yields independent events. It is therefore not surprising that the bounds obtained by using this repeated compression method are in close connection to those obtained by the local lemma. 

\section{Application to other Ramsey objects}
{

Naturally the same method can be applied to establish other Ramsey-theoretic lower bounds for example for Ramsey numbers or transitive subtournaments. In both cases the factor $e/4$ vanishes asymptotically as $n$ appears superlinear in the equations. The bounds therefore match (up to the $k$-th root of a constant) the bounds obtained with the Lov\'{a}sz local lemma~\cite{alon_spencer_92}.

}
	
}

\bibliography{the_compression_method}

\bibliographystyle{abbrv}

\end{document}